\documentclass[12pt]{amsart}
\usepackage[latin1]{inputenc}
\usepackage{amscd}
\usepackage{latexsym}
\usepackage{pstcol,pst-plot,pst-3d}
\usepackage{mathptmx}
\usepackage{multicol}

\psset{unit=0.7cm,linewidth=0.8pt,arrowsize=2.5pt 4}
\def\vertex{\pscircle[fillstyle=solid,fillcolor=black]{0.03}}
\newpsstyle{fatline}{linewidth=1.5pt}
\newpsstyle{fyp}{fillstyle=solid,fillcolor=verylight}
\definecolor{verylight}{gray}{0.97} 
\definecolor{light}{gray}{0.9}
\definecolor{medium}{gray}{0.85}

\def\opn#1#2{\def#1{\operatorname{#2}}} 
\opn\chara{char}
\opn\length{\ell}
\opn\pd{pd}
\opn\rk{rk}
\opn\projdim{proj\,dim}
\opn\injdim{inj\,dim}
\opn\rank{rank}
\opn\depth{depth}
\opn\grade{grade}
\opn\height{height}
\opn\embdim{emb\,dim}
\opn\codim{codim}

\opn\Tr{Tr}
\opn\bigrank{big\,rank}
\opn\superheight{superheight}\opn\lcm{lcm}
\opn\trdeg{tr\,deg}%
\opn\reg{reg}
\opn\lreg{lreg}
\opn\skel{skel}
\opn\div{div}
\opn\Div{Div}
\opn\cl{cl}
\opn\Cl{Cl}
\opn\Spec{Spec}
\opn\Supp{Supp}
\opn\supp{supp}
\opn\Sing{Sing}
\opn\Ass{Ass}
\opn\Ann{Ann}
\opn\Rad{Rad}
\opn\Soc{Soc}
\opn\Ker{Ker}
\opn\Coker{Coker}
\opn\Im{Im}
\opn\Hom{Hom}
\opn\Tor{Tor}
\opn\Ext{Ext}
\opn\End{End}
\opn\Aut{Aut}
\opn\id{id}

\opn\nat{nat}
\opn\pff{pf}
\opn\Pf{Pf}
\opn\GL{GL}
\opn\SL{SL}
\opn\mod{mod}
\opn\ord{ord}
\opn\aff{aff}
\opn\con{conv}
\opn\relint{relint}
\opn\st{st}
\opn\lk{lk}
\opn\cn{cn}
\opn\core{core}
\opn\vol{vol}
\opn\link{link}
\opn\star{star}
\opn\skel{skel}
\opn\gr{gr}

\def\pot#1#2{#1[\kern-0.28ex[#2]\kern-0.28ex]}

\opn\dirlim{\underrightarrow{\lim}}
\opn\inivlim{\underleftarrow{\lim}}
\let\sect=\cap

\let\Dirsum=\bigoplus

\let\to=\rightarrow
\let\To=\longrightarrow
\def\Implies{\ifmmode\Longrightarrow \else
     \unskip${}\Longrightarrow{}$\ignorespaces\fi}
\def\implies{\ifmmode\Rightarrow \else
     \unskip${}\Rightarrow{}$\ignorespaces\fi}
\def\iff{\ifmmode\Longleftrightarrow \else
     \unskip${}\Longleftrightarrow{}$\ignorespaces\fi}

\let\:=\colon
\newtheorem{Theorem}{Theorem}[section]
\newtheorem{Lemma}[Theorem]{Lemma}
\newtheorem{Corollary}[Theorem]{Corollary}
\newtheorem{Proposition}[Theorem]{Proposition}

\let\epsilon\varepsilon
\let\phi=\varphi
\let\kappa=\varkappa
\def\qed{\ifhmode\textqed\fi
   \ifmmode\ifinner\quad\qedsymbol\else\dispqed\fi\fi}
\def\textqed{\unskip\nobreak\penalty50
    \hskip2em\hbox{}\nobreak\hfil\qedsymbol
    \parfillskip=0pt \finalhyphendemerits=0}
\def\dispqed{\rlap{\qquad\qedsymbol}}

\opn\inii{in}
\opn\inim{inm}
\opn\rate{rate}
\begin{document}
\title{Dirac's theorem on chordal graphs and Alexander duality }
\author{J\"urgen Herzog, Takayuki Hibi and Xinxian Zheng}
\address{J\"urgen Herzog, Fachbereich Mathematik und
Informatik,  
Universit\"at Duisburg-Essen, 45117 Essen, Germany}
\email{juergen.herzog@uni-essen.de}
\address{Takayuki Hibi, Department of Pure and Applied Mathematics, Graduate School 
of Information Science and Technology,
Osaka University, Toyonaka, Osaka 560-0043, Japan}
\email{hibi@math.sci.osaka-u.ac.jp}
\address{Xinxian Zheng, Fachbereich Mathematik und
Informatik,  
Universit\"at Duisburg-Essen, 45117 Essen, Germany}
\email{xinxian.zheng@uni-essen.de}
\date{}
\maketitle
\begin{abstract}
By using Alexander duality on simplicial complexes we give a new and algebraic proof of Dirac's theorem on chordal 
graphs. 
\end{abstract}

\section*{Introduction}
One of the fascinating results in classical graph theory is Dirac's theorem \cite{D} 
on chordal graphs.  Recall that a finite graph $G$ is chordal if each cycle of $G$ 
of length $\geq 4$ has a chord.  
Using our terminology `quasi-trees' (introduced in the beginning of Section $2$), 
Dirac proved that a finite graph $G$ is chordal 
if and only if $G$ is the $1$-skeleton of a quasi-tree.

In commutative algebra, the chordal graph first appeared in the work of Fr\"oberg \cite{F}.  Let $S = K[x_1, \ldots, 
x_n]$ denote the polynomial ring in $n$ variables
over a field $K$.  Given a finite graph $G$ on 
$[n] = \{ 1, \ldots, n \}$, we associate the monomial ideal $I(G) \subset S$,
called the edge ideal of $G$, 
generated by those monomials $x_ix_j$ such that $\{ i, j \}$ is an edge of $G$.
In \cite{F} it is proved that $I(G)$ has a linear resolution if and only if the complementary graph $\bar{G}$ of $G$ 
is chordal. 
Recently, in \cite{HHZ} it is proved that if $\bar{G}$ is chordal, then all powers of 
$I(G)$ have linear resolutions.  The Dirac's theorem plays an essential role in \cite{HHZ}.

Explaining the results of \cite{HHZ} to David Eisenbud when the first two authors visited MSRI, he expressed his 
opinion that the quasi-trees appearing in Fr\"oberg's theorem on edge ideals should be related via Alexander duality 
to trees that are naturally attached to the relation matrix of a monomial ideal which is perfect of codimension 2. 
The main purpose of this paper is to show that this is indeed the case, and thereby giving a new and algebraic proof 
of Dirac's theorem.

The present paper is organized as follows.  In Section $1$ we discuss basic concepts related to simplicial complexes 
such as Stanley--Reisner ideals, facet ideals, Alexander duality, skeletons and flag complexes.

The crucial Lemma \ref{relationmatrix} is proved in Section 2, where it is shown that a quasi-tree is characterized in 
terms of the Taylor relations of a certain monomial ideal. Combining this fact with the Hilbert--Burch theorem we show 
in Corollary \ref{projdim} that a simplicial complex $\Delta$ is a quasi-tree if and only if  the facet ideal 
$I(\Delta^c)$ has  projective dimension one, where $\Delta^c$ is the simplicial complex whose facets are the 
complements of the facets of $\Delta$. 

Our algebraic proof of Dirac's theorem is presented in Section 3, see Theorem \ref{Dirac}. There we also discuss a 
sort of higher Dirac theorem. Finally in Section 4 we extend the main result of \cite{HHZ} showing that all powers of 
non-skeleton facet ideals of a quasi-tree have  linear resolutions.

\section{Stanley--Reisner ideals and facet ideals}
Let $S = K[x_1, \ldots, x_n]$ denote the polynomial ring in $n$ variables
over a field $K$.
Write $[n]$ for the finite set $\{ 1, \ldots, n \}$ and
${[n] \choose i}$ the set of all $i$-element subsets of $[n]$.

Let $\Delta$ be a simplicial complex on the vertex set $[n]$. Thus $\Delta$ is a collection of subsets of $[n]$ such 
that (i) $\{ i \} \in \Delta$ for all $i \in [n]$ and
(ii) if $F \in \Delta$ and $G \subset F$,
then $G \in \Delta$.  Each element $F \in \Delta$
is called a {\em face} of $\Delta$.  
The dimension of a face $F$ is $|F| - 1$.
Here $|F|$ is the cardinality of a finite set $F$.
The dimension of $\Delta$ is $\dim \Delta = \max\{|F| \: F \in \Delta \} - 1$.
A {\em facet} of $\Delta$ is a maximal face of $\Delta$.  
A {\em nonface} of $\Delta$ is a subset $G$ of $[n]$
with $G \not\in \Delta$.
Let ${\mathcal F}(\Delta)$ denote the set of facets of $\Delta$. 
A simplicial complex $\Delta$ is called {\em pure} if
all the facets of $\Delta$ have the same cardinality.

Naturally attached to $\Delta$ are two squarefree monomial ideals in $S$. The first, more commonly known ideal, is the 
{\em Stanley--Reisner ideal} $I_{\Delta}$, which is generated by all monomials $x_F$ such that $F\not\in \Delta$. Here 
$x_F=x_{i_1}\cdots x_{i_k}$ for $F=\{i_1,\ldots, i_k\}$. The second is the so-called {\em facet ideal} $I(\Delta)$ 
which  is generated by all monomials $x_F$ where $F$ is a facet of $\Delta$.  In case $\Delta=G$ is a graph, $I(G)$ is 
called the edge ideal of $G$. Suppose ${\mathcal F}(\Delta)=\{F_1,\ldots, F_m\}$. Then we write $\Delta=\langle 
F_1,\ldots,F_m\rangle$, and we have $I(\Delta)=(x_{F_1},\ldots,x_{F_m})$.

In this section we want to discuss the relationship between these two ideals.

Suppose $\Delta$ is a pure $(d-1)$-dimensional simplicial complex.  We then define
\[
\bar{\Delta}=\langle F\: F\not\in \Delta, |F|\in {[n] \choose d} \rangle. 
\]

Recall that that the $i$-skeleton of a simplicial complex $\Delta$ is the simplicial complex 
$\skel_{\Delta}(i)$ whose facets are the $i$-dimensional  faces of $\Delta$. 

We have the following very simple

\begin{Lemma}
\label{trivial}
Let $\Delta$ be a $(d-1)$-dimensional pure simplicial complex, and let $\Gamma$ be the simplicial complex such that 
$I(\Delta)=I_\Gamma$. Then 
\[
\bar{\Delta}=\skel_\Gamma(d-1). 
\]
\end{Lemma}

\begin{proof}
Let $F\in{\mathcal F}(\bar{\Delta})$, then $F\not\in\Delta$. Therefore $x_F\not\in I(\Delta)$, and hence $x_F\not\in 
I_\Gamma$. This means that $F\in \Gamma$. Since $|F|=d$, this implies that 
$F\in \skel_\Gamma(d-1)$. The converse inclusion is proved similarly.
\end{proof}

Next we express the Stanley--Reisner ideal of the Alexander dual of a simplicial complex $\Delta$ in terms of a facet 
ideal. Recall that the simplicial complex 
\[
\Delta^\vee=\{[n]\setminus F\: F\not\in\Delta\}
\]
is called the {\em Alexander dual}  of $\Delta$. It is easy to see that $(\Delta^\vee)^\vee=\Delta$.  

We also define
\[
\Delta^c= \langle [n]\setminus F\: F\in {\mathcal F}(\Delta)\rangle. 
\]
We denote $[n]\setminus F$ by $F^c$. As usual, we use $G(I)$ to denote the unique minimal generating system of
the monomial ideal $I$.

\begin{Lemma}
\label{true}
Let $\Delta$ be a simplicial complex. Then 
\[
I_{\Delta^\vee}=I(\Delta^c).
\]
\end{Lemma}

\begin{proof}
By definition, 
$\Delta^\vee=\langle F^c \: F \text { is a minimal nonface of }\Delta \rangle$.  
Furthermore,  $x_G\in G(I_{\Delta^\vee})$ if and only if $G$ is a minimal subset of $[n]$ such that 
$G\not\in\Delta^\vee$. This means that  $G^c$ does not contain any minimal nonface of $\Delta$, and for any proper 
subset $H$ of $G$, the complement $H^c$ contains a minimal nonface of $\Delta$. This is equivalent to say that any 
subset of $G^c$ is a face of $\Delta$, and for any proper subset $H$ of $G$, 
$H^c$ is not a face of $\Delta$. In another words, $G^c$ is a facet of $\Delta$. Hence 
$I(\Delta^c)=I_{\Delta^\vee}$, as required.
\end{proof}

Recall that a simplicial complex is called {\em flag}, if all minimal nonfaces consist of two elements, equivalently, 
$I_\Delta$ is generated by quadratic monomials.  We also consider the simplex on $[n]$ as a flag complex. Note that if 
$\Delta$ has only two facets, then $\Delta$ is flag.

\begin{Proposition}
\label{skeleton}
Let $\Sigma$ be a flag complex with $n$ vertices, and let $\Delta$ and $\Delta'$ be the simplicial complexes defined 
by 
\[
I_\Delta=I(\overline{\skel_\Sigma(\ell)})\quad\text{and}\quad 
I_{\Delta'}=I(\overline{\skel_\Sigma(1)}).
\]
Then $\Delta^\vee=\skel_{(\Delta')^\vee}(n-\ell-2)$.
\end{Proposition}

\begin{proof} By Lemma \ref{true} we have $\Delta^\vee=(\overline{\skel_\Sigma(\ell)})^c$ and 
$(\Delta')^\vee=(\overline{\skel_\Sigma(1)})^c$. Since $\Sigma$ is flag, any facet of  $\overline{\skel_\Sigma(\ell)}$ 
contains a nonedge of $\Sigma$ which is a facet of 
$\overline{\skel_\Sigma(1)}$. Therefore, any facet of $\Delta^\vee$ is a face of $(\Delta')^\vee$. It is clear that 
the facets of $\Delta^\vee$ are all of dimension $n-\ell-2$, so that 
$\Delta^\vee\subset\skel_{(\Delta')^\vee}(n-\ell-2)$.

On the other hand, for any $(n-\ell-2)$-dimensional face $F$ of $(\Delta')^\vee$ its complementary set $F^c$ contains 
one nonedge of $\Sigma$. Therefore, $F^c\in \overline{\skel_\Sigma(\ell)}$ and hence $F$ is a facet of $\Delta^\vee$. 
\end{proof}

We quote the following two results relating combinatorial or algebraic properties of  a simplicial complex $\Delta$ to 
algebraic properties of the Alexander dual of 
$\Delta^\vee$.

\begin{Theorem}
\label{duality}
Let $K$ be field, and $\Delta$ be a simplicial complex. Then 
\begin{enumerate}
\item[(a)] {\em (Eagon--Reiner \cite{ER})} $K[\Delta]$ is  Cohen-Macaulay \iff
 $I_{\Delta^\vee}$ has a linear resolution.
\item[(b)] {\em (Terai \cite{T})} $\projdim K[\Delta]=\reg(I_{\Delta^\vee})$.
\item[(c)]  $\Delta$ is shellable \iff 
$I_{\Delta^\vee}$ has linear quotients.
\end{enumerate}
\end{Theorem}

For the convenience of the reader we give the easy proof of  statement (c): recall that $\Delta$ is called {\em 
shellable} if $\Delta$ is pure and there is an order $F_1,\ldots,F_m$ of the facets of $\Delta$  (called a {\em 
shelling order}), such that for all  $0<j<i$ and $x\in F_i\setminus F_j$, there exists $k<i$ such that $F_i\setminus 
F_k=\{x\}$, while 
an ideal $I$ is said to have {\em linear quotients}, if $I=(f_1,\ldots, f_m)$ and for all $i>0$  the colon ideals 
$(f_1,\ldots,f_{i-1}):f_i$ are generated by linear forms.

For a monomial ideal $I$ we require that the $f_i$ belong to the unique minimal set of monomial generators $G(I)$ of 
$I$. Then $I$ has linear quotients if for all $i>1$, and any $j<i$, there exists $k<i$ such that $f_k/[f_i,f_k]$ is a 
monomial of degree 1, say $x_\ell$,  and $x_\ell|f_j$.
Here $[f_i,f_k]$ denotes the greatest common divisor of $f_i$ and $f_k$.
By Lemma \ref{true} one has $I_{\Delta^\vee}=(x_{F_1^c},\ldots, x_{F_m^c})$,  the equivalence of the statements in (c) 
are obvious.

It is well known that $K[\Delta]$ is Cohen-Macaulay for any field $K$, if $\Delta$ is shellable, see for instance 
\cite{BH}, and it is  easy to see that an ideal with linear quotients has a linear resolution.

\begin{Corollary}
\label{application}
 Let $\Sigma$ be a flag complex, and let $\Delta$ and $\Delta'$ be the simplicial complexes defined by 
$I_\Delta=I(\overline{\skel_\Sigma(\ell)})$ and $I_{\Delta'}=I(\overline{\skel_\Sigma(1)})$. Suppose that 
$I_{\Delta'}$ has linear quotients, then so does $I_\Delta$.
\end{Corollary}

\begin{proof}
It follows from Theorem \ref{duality}(c) that $(\Delta')^\vee$ is shellable. Since $\Delta^\vee$ is a skeleton of 
$(\Delta')^\vee$, the next lemma implies that $\Delta^\vee$ is shellable, too. Applying again Theorem 
\ref{duality}(c), the assertion follows.
\end{proof}

\begin{Lemma}
\label{shellable}
Let $\Delta$ be a shellable complex with $\dim \Delta=d-1$. Then for each $1 \leq i < d$ the $i$-skeleton 
$\skel_\Delta(i)$ of $\Delta$ is shellable.
\end{Lemma}

\begin{proof}
Let $i<d-1$.
Fix a shelling $F_1, \ldots, F_m$ of the facets of $\Delta$.
If $m=1$, i.e., $\Delta$ is the simplex on $[n]$, then ${\mathcal F}(\skel_\Delta(i)) = {[n] \choose i + 1}$, and 
$\skel_\Delta(i)$ is shellable.  
Let $m>1$ and $\Delta'=\Delta\setminus \{F_m\}$.  By using induction on $m$, we may assume that 
$\skel_{\Delta'}(i)$ is shellable.  Let $V\subset [n]$ denote the set of those $b\in [n]$ such that there is 
$1\leq s < m$ with $\dim(F_s \cap F_m)=d-2$ and with $F_m \setminus F_s = \{b\}$.  
It then follows that a subset $G \in {[n] \choose i + 1}$ belongs to 
${\mathcal F}(\skel_{\Delta}(i))\setminus {\mathcal F}(\skel_{\Delta'}(i))$ if and only if 
$V\subset G\subset F_m$. Hence the simplicial complex $\Gamma$ with 
${\mathcal F}(\Gamma)={\mathcal F}(\skel_{\Delta}(i)) \setminus {\mathcal F}(\skel_{\Delta'}(i))$
turns out to be shellable.

Let $G_1, G_2, \ldots, G_p$ be a shelling of the facets of $\skel_{\Delta'}(i)$ and
$G_{p+1}, \ldots, G_{q}$ a shelling of $\Gamma$. We claim that 
$G_1, G_2, \ldots, G_p, G_{p+1}, \ldots, G_{q}$ is a shelling of $\skel_\Delta(i)$.
In fact, let $1\leq j\leq p < k \leq q$ and $G_j \subset F_s$ with $s<m$.
Then there is $s'< m$ with $\dim(F_{s'} \cap F_m)=d-2$ such that
$F_s\cap F_m\subset F_{s'}\cap F_m$. Let $F_{s'} \setminus F_m =\{a\}$ and
$F_m \setminus F_{s'}=\{b\}$. Since $p<k$, one has $b\in G_k$.
Let $G_{k'}=(G_k \setminus\{b\})\cup\{a\}$ with $k'\leq p$.
Then $G_{k'} \cap G_k=G_k \setminus \{ b \} \in {[n] \choose i}$.
Since $b \not\in F_s$, one has $b \not\in G_j$.
Hence $G_j \cap G_k\subset G_{k'} \cap G_k$, as desired.
\end{proof}

\section{Quasi-trees and relation trees of ideals of projective dimension 1}

Let $\Delta$ be a simplicial complex. A facet $F\in {\mathcal F}(\Delta)$  is called a {\em leaf}, if either $F$ is 
the only facet of $\Delta$, or there exists  $G\in {\mathcal F}(\Delta)$, $G\neq F$ such that $H\sect F\subset G\sect 
F$ for each $H\in{\mathcal F}(\Delta)$ with $H\neq F$. A 
facet $G$ with this property is called a {\em branch} of $F$. A vertex $i$  of $\Delta$ is called a {\em  free vertex} 
if $i$ belongs to precisely one facet. 

Faridi \cite{Fa} calls $\Delta$ a tree if each simplicial complex generated by a subset of the facets of $\Delta$ has 
a leaf, and Zheng \cite{Z} calls $\Delta$ a {\em quasi-tree} if there exists a labeling $F_1,\ldots, F_m$ of the 
facets such that for all $i$ the facet $F_i$ is a leaf of the subcomplex $\langle F_1,\ldots, F_i\rangle$. We call 
such a labeling a {\em leaf order}. It is obvious that any tree is a quasi-tree, but the converse is not true. For us 
however the quasi-trees are important.

Let $\Delta$ be a simplicial complex on $[n]$ with
${\mathcal F}(\Delta) = \{ F_1, \ldots, F_t \}$.
We introduce the ${t \choose 2} \times t$ matrix
\[
M_{\Delta} = ( a^{(i,j)}_{k} )_
{1 \leq i < j \leq t, \, 1 \leq k \leq t}
\]
whose entries $a^{(i,j)}_{k} \in S$ are
$a^{(i,j)}_{i} = x_{F_i \setminus F_j}$, 
$a^{(i,j)}_{j} = x_{F_j \setminus F_i}$, and
$a^{(i,j)}_{k} = 0$
if $k \not\in \{ i, j \}$
for all $1 \leq i < j \leq t$ and for all $1 \leq k \leq t$.

\begin{Lemma} 
\label{relationmatrix}
A simplicial complex $\Delta=\langle F_1,\ldots,F_t \rangle$ on $[n]$ 
is a quasi-tree if and only if the matrix $M_{\Delta}$ contains a $(t - 1) \times t$ submatrix $M^\sharp_{\Delta}$ 
with the property that, for each $1 \leq j \leq t$, 
if $M^\sharp_{\Delta}(j)$
is the $(t - 1) \times (t - 1)$ submatrix of $M^\sharp_{\Delta}$
obtained by removing the $j$th column from $M^\sharp_{\Delta}$,
then $|\det(M^\sharp_{\Delta}(j))| = x_{[n]} / x_{F_j}$. 
\end{Lemma}

\begin{proof}
({\bf ``only if''})
Let $\Delta$ be a quasi-tree on $[n]$ and fix a leaf ordering
$F_1, \ldots, F_t$ of the facets of $\Delta$.
Let $t > 1$.
Let $F_k$ with $k \neq t$
be a branch of $F_t$
and $\Delta' = \Delta \setminus F_t$.
Since $\Delta'$ is a quasi-tree, by assumption of induction,
it follows that   
$M_{\Delta}$ contains a $(t - 2) \times t$ submatrix
$M'$ with the property that, for each $1\leq j<t$, 
if $M'(j,t)$
is the $(t - 2) \times (t - 2)$ submatrix of $M'$
obtained by removing the $j$st and $t$th columns from $M'$,
then $|\det(M'(j,t))| = 
x_{[n] \setminus (F_t \setminus F_k)} / x_{F_j}$.
Let $M^\sharp_{\Delta}$ denote the 
$(t - 1) \times t$ submatrix of $M_{\Delta}$
obtained by adding the $(k,t)$th row to $M'$.
Since $a^{(k,t)}_{t} = x_{F_t \setminus F_k}$,
it follows that, 
for each $1\leq j<t$, one has 
$|\det(M^\sharp_{\Delta}(j))| = x_{[n]} / x_{F_j}$.
Moreover, since
$|\det(M^\sharp_{\Delta}(t))|
= |x_{F_k \setminus F_t} \det(M'(k,t))|$,
one has 
$|\det(M^\sharp_{\Delta}(t))|
= x_{[n]} / x_{F_t}$. 

({\bf ``if''})
Now, suppose that the matrix $M_{\Delta}$ contains 
a $(t - 1) \times t$ submatrix $M^\sharp_{\Delta}$ 
with the property that, for each $1 \leq j \leq s$, 
if $M^\sharp_{\Delta}(j)$ is the $(t - 1) \times (t - 1)$ submatrix 
of $M^\sharp_{\Delta}$ obtained by removing $j$th column from 
$M^\sharp_{\Delta}$,
then $|\det(M^\sharp_{\Delta}(j))| = x_{[n]} / x_{F_j}$.
Let $\Omega$ denote the subgraph on $[t]$
whose edges are those $\{ i, j \}$ with $1 \leq i < j \leq t$
such that the $(i, j)$th row of $M_{\Delta}$ belongs to
$M^\sharp_{\Delta}$.    
Then $\Omega$ contains no cycles.
To see why this is true, if $C$ is a cycle of $\Omega$
with $E(C)$ its edge set.
If $\{ i_0, j_0 \} \in E(C)$,
then in the matrix  
$M^\sharp_{\Delta}(i_0)$,
the $(i, j)$th rows with $\{i, j\} \in E(G)$ 
are linearly dependent.
Thus $\det(M^\sharp_{\Delta}(i_0)) = 0$.
This is impossible.
Hence $\Omega$ contains no cycles.
Since the number of edges of $\Omega$ is $t - 1$,
it follows that $\Omega$ is a tree, i.e., a connected graph
without cycles.  Hence there is a column 
of $M^\sharp_{\Delta}$ which contains exactly one nonzero entry.
Suppose, say, that the $t$st column contains exactly one nonzero 
entry and the $(k,t)$th row of $M_{\Delta}$ appears in 
$M^\sharp_{\Delta}$.  Then, for each $1\leq j<t$, 
the monomial $x_{F_t \setminus F_k}$ divides
$|\det(M^\sharp_{\Delta}(j))|$.  Hence
$(F_t \setminus F_k) \cap F_j = \emptyset$
for all $1\leq j<t$.
It then follows that $F_t$ is a leaf of $\Delta$
and $F_k$ is a branch of $F_t$.
Let $\Delta' = \Delta \setminus F_t$
and $M^\sharp_{\Delta'}$ the $(t-2) \times (t-1)$
submatrix of $M_{\Delta'}$ which is obtained by removing
the $(k,t)$th row and the $t$st column from 
$M^\sharp_{\Delta}$. 
Since $\Delta'$ is a simplicial complex on 
$[n]\setminus (F_t\setminus F_k)$
and since $x_{F_t \setminus F_k}
(x_{[n] \setminus (F_t\setminus F_k)} / x_{F_j})
= x_{[n]} / x_{F_j}$ for each $1\leq j<t$,
working with induction on $t$, it follows that
$\Delta'$ is a quasi-tree.  Hence $\Delta$ is a quasi-tree.
\end{proof}

Let $I$ be an arbitrary monomial ideal with $G(I)=\{u_1,\ldots,u_t\}$, and let $T$ be the Taylor complex associated 
with $I$. Then $T_i=S^{\binom{t}{i}}$, and the matrix $A_I$ representing the differential $T_2\to T_1$ is a 
$\binom{t}{2}\times t$-matrix. To be more precise, if $T_1=\Dirsum_{i=1}^tSe_i$, then $T_2=\Dirsum_{i<j}Se_i\wedge 
e_j$,  and $\partial(e_i\wedge e_j)=
u_{ji}e_i-u_{ij}e_j$, where $u_{ij}=u_i/[u_i,u_j]$ for all $i,j\in[t]$  with $i\neq j$.

Note that for any simplicial  complex $\Delta$ we have $M_{\Delta}= A_{I(\Delta^c)}$, because if $u_i=x_{F^c_i}$ and 
$u_j=x_{F^c_j}$, then $u_{ji}=x_{F^c_j\setminus F^c_i}=x_{F_i\setminus F_j}$.

Assume now that $I$ has projective dimension 1, and that the elements of $G(I)$   have no common factor. Then $I$ is 
perfect of codimension 2. 

A subset $R$ of the Taylor relations is called {\em irreducible} if $R$ generates the first syzygy module 
${\text{syz}}_1(I)$ of $I$, but no proper subset of  $R$ generates ${\text{syz}}_1(I)$. Fortunately it is known (see
\cite[Corollary 5.2]{BH1}) that an irreducible subset of the Taylor relations is in fact a minimal system of 
generators of ${\text{syz}}_1(I)$. In particular it follows that we can always choose a minimal free resolution 
\begin{eqnarray*}
\begin{CD}
0@>>> S^{t-1}@>\phi >> S^t@>>> I@>>> 0
\end{CD}
\end{eqnarray*}
such that the rows of the matrix of $\phi$ correspond to Taylor relations. However the choice of
an irreducible set $R$ of  Taylor relations  is in general not unique.

For example, let $I=(x_4x_5x_6, x_1x_5x_6, x_1x_2x_6, x_1x_2x_5)$. Then $\phi$ can be represented by the matrix 
\[
\begin{pmatrix}
x_1 & -x_4&0&0\\
0&x_2&-x_5&0\\
0&x_2&0&-x_6
\end{pmatrix},
\]
or by 
\begin{center}
$
\begin{pmatrix}
x_1 & -x_4&0&0\\
0&x_2&-x_5&0\\
0&0&x_5&-x_6
\end{pmatrix}
$
\ \ \ or  \ \ \ 
$
\begin{pmatrix}
x_1 & -x_4&0&0\\
0&x_2&0&-x_6\\
0&0&x_5&-x_6
\end{pmatrix}.
$
\end{center}
\medskip
Nevertheless for a given choice $R$ of $t-1$ Taylor relations which generate ${\text{syz}}_1(I)$ we can define a 
($1$-dimensional) tree $\Omega$ as in the proof of \ref{relationmatrix} with  
$\{i,j\}\in E(\Omega)$  if $u_{ji}e_i-u_{ij}e_j\in R$ for $i<j$. We call $\Omega$ the {\em relation tree} of $R$. This 
relation tree was first considered in \cite[Remark 6.3]{BH1}.

In the above example the relation tree for the first matrix is
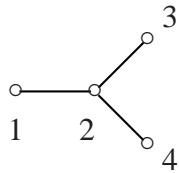
\begin{center}

\psset{unit=1.0cm}
\begin{pspicture}(0,0)(3,3)
 \psline(0.5,1.5)(1.45,1.5)
 \psline(1.55,1.45)(2.15,0.85)
 \psline(1.55,1.55)(2.15,2.15)
 \rput(0.45,1.5){$\circ$}
 \rput(1.5,1.5){$\circ$}
 \rput(2.2,0.8){$\circ$}
 \rput(2.2,2.2){$\circ$}
 \rput(0.45,1){1}
 \rput(1.4,1){2}
 \rput(2.5,2.5){3}
 \rput(2.5,0.6){4}
\end{pspicture}

\end{center}

while for the other matrices it is 

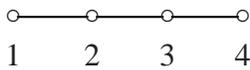
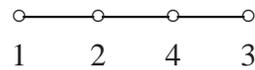
\begin{center}
\psset{unit=1.0cm}
\begin{pspicture}(0,0)(12,2.5)
 \psline(0.5,1.5)(1.45,1.5)
 \psline(1.55,1.5)(2.45,1.5)
 \psline(2.55,1.5)(3.45,1.5)
 \rput(0.45,1.5){$\circ$}
 \rput(1.5,1.5){$\circ$}
 \rput(2.5,1.5){$\circ$}
 \rput(3.5,1.5){$\circ$}
 \rput(0.45,1){1}
 \rput(1.5,1){2}
 \rput(2.5,1){3}
 \rput(3.5,1){4}
 \rput(6,1.5){$\text{or}$}
 \psline(8.5,1.5)(9.45,1.5)
 \psline(9.55,1.5)(10.45,1.5)
 \psline(10.55,1.5)(11.45,1.5)
 \rput(8.45,1.5){$\circ$}
 \rput(9.5,1.5){$\circ$}
 \rput(10.5,1.5){$\circ$}
 \rput(11.5,1.5){$\circ$}
 \rput(8.45,1){1}
 \rput(9.5,1){2}
 \rput(10.5,1){4}
 \rput(11.5,1){3}
\end{pspicture}
\end{center}

Next we want to describe how the generators $u_i$ of $I$ can be computed from the $u_{ij}$ and the relation trees. To 
this end we introduce for each $i=1,\ldots, t$ an orientation to make $\Omega$ a directed graph which  we denote 
$\Omega_i$. We fix some vertex $i$. Let $j$ be any other vertex of $\Omega$. Since $\Omega$ is a tree there is a 
unique directed walk from $i$ to $j$. This defines the orientation of the edges along this walk. The following picture 
explains this for the first of our relation trees in the above example.

\psset{unit=1.0cm}
\begin{pspicture}(0,0)(3,3)
 \psline{->}(0.5,1.5)(1.45,1.5)
 \psline{->}(1.55,1.45)(2.15,0.85)
 \psline{->}(1.55,1.55)(2.15,2.15)
 \rput(0.45,1.5){$\circ$}
 \rput(1.5,1.5){$\circ$}
 \rput(2.2,0.8){$\circ$}
 \rput(2.2,2.2){$\circ$}
 \rput(0.45,1){1}
 \rput(1.4,1){2}
 \rput(2.5,2.5){3}
 \rput(2.5,0.6){4}

 \psline{<-}(3.5,1.5)(4.45,1.5)
 \psline{->}(4.55,1.45)(5.15,0.85)
 \psline{->}(4.55,1.55)(5.15,2.15)
 \rput(3.45,1.5){$\circ$}
 \rput(4.5,1.5){$\circ$}
 \rput(5.2,0.8){$\circ$}
 \rput(5.2,2.2){$\circ$}
 \rput(3.45,1){1}
 \rput(4.4,1){2}
 \rput(5.5,2.5){3}
 \rput(5.5,0.6){4}

 \psline{<-}(6.5,1.5)(7.45,1.5)
 \psline{->}(7.55,1.45)(8.15,0.85)
 \psline{<-}(7.55,1.55)(8.15,2.15)
 \rput(6.45,1.5){$\circ$}
 \rput(7.5,1.5){$\circ$}
 \rput(8.2,0.8){$\circ$}
 \rput(8.2,2.2){$\circ$}
 \rput(6.45,1){1}
 \rput(7.4,1){2}
 \rput(8.5,2.5){3}
 \rput(8.5,0.6){4}
 
 \psline{<-}(9.5,1.5)(10.45,1.5)
 \psline{<-}(10.55,1.45)(11.15,0.85)
 \psline{->}(10.55,1.55)(11.15,2.15)
 \rput(9.45,1.5){$\circ$}
 \rput(10.5,1.5){$\circ$}
 \rput(11.2,0.8){$\circ$}
 \rput(11.2,2.2){$\circ$}
 \rput(9.45,1){1}
 \rput(10.4,1){2}
 \rput(11.5,2.5){3}
 \rput(11.5,0.6){4}
\end{pspicture}

By  the Hilbert--Burch theorem  \cite[Theorem 1.4.17]{BH} one has
\[
u_i=(-1)^{i}\det(A_i)\quad\text{for}\quad i=1,\ldots, t,
\]
where the matrix $A_i$ is obtained from the relation matrix $A$ of $I$ by deleting the $i$th column of $A$. Computing 
$\det(A_i)$ by the determinantal expansion formula as in the proof of Lemma \ref{relationmatrix} one sees that 
\[
u_i=\prod_{(k,j)}u_{kj},
\]
where the product is taken over all oriented edges $(k,j)$ of $\Omega_i$.

\begin{Corollary}
\label{projdim}
A simplicial complex $\Delta$  is a quasi-tree if and only if 
 $\projdim I(\Delta^c) = 1$.
\end{Corollary}

\begin{proof} 
Let ${\mathcal F}(\Delta)=\{F_1,\ldots, F_t\}$. By Lemma \ref{relation matrix},
the simplicial complex $\Delta$ is a quasi-tree if and only if $M_{\Delta}$ contains a $(t-1)\times t$ submatrix  
$M_{\Delta}^\#$ whose ideal of maximal minors is $I(\Delta^c)$. 
Hence, if $\Delta$ is a quasi-tree, the Hilbert--Burch theorem  implies that 
$\projdim I(\Delta^c)=1$. Conversely, suppose $\projdim I(\Delta^c)=1$, and let $A$ be a $(t-1)\times t$ relation 
matrix of this ideal consisting of Taylor relations. By the Hilbert-Burch theorem, $I(\Delta^c)$ is the ideal of 
maximal minors of $A$. Since $M_{\Delta}=M_{\Delta^c}$, it follows that $A$ is a submatrix of $M_\Delta$.  Hence 
$\Delta$ is a quasi-tree.
\end{proof}

In our example $I$ may be viewed as $I=I(\Delta^c)$ where the facets of $\Delta$ are 
\[
\{\{a,b,c\},\{b,c,d\},\{c,d,e\},\{c,d,f\}\}.
\] 
See the following picture:

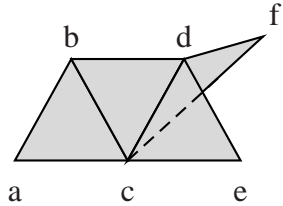
\begin{center}
\psset{unit=1.5cm}
\begin{pspicture}(0,0)(3,2)
\pspolygon[style=fyp, fillcolor=medium](0.5,0.5)(1.5,0.5)(1,1.4)
 \pspolygon[style=fyp, fillcolor=medium](1.5,0.5)(1,1.4)(2,1.4)
 \pspolygon[style=fyp, fillcolor=medium](1.5,0.5)(2,1.4)(2.7,1.6)
  \pspolygon[style=fyp, fillcolor=medium](1.5,0.5)(2,1.4)(2.5,0.5)
 \psline[linestyle=dashed](1.5,0.5)(2.3,1.23)
 \rput(0.5,0.2){a}
 \rput(1,1.6){b}
 \rput(1.5,0.2){c}
 \rput(2,1.6){d}
 \rput(2.5,0.2){e}
 \rput(2.8,1.8){f}
\end{pspicture}
\end{center}
This is a quasi-tree, as it should be by Corollary \ref{projdim}.

Inspecting the proof of Lemma \ref{relationmatrix}, we see that  all possible relation trees $\Omega$ of $I(\Delta^c)$ 
can be recovered from the quasi-tree  $\Delta=\langle F_1,\ldots, F_m\rangle$ as follows: start with some leaf  $F_i$  
of $\Delta$,  and let $F_j$ be a branch of $F_i$.  Then $\{i,j\}$ will be an edge of $\Omega$. According to Corollary 
\ref{hibifound}, $\langle {\mathcal F}(\Delta)\setminus \{F_i\}\rangle$ is again a quasi-tree. Then remove the leaf 
$F_i$, and continue in the same way with the remaining quasi-tree in order to find the other edges of $\Omega$. Of 
course, at each step of the procedure there may be different choices. This gives us the different possible relation 
trees. 

Geometrically a relation tree is obtained from a given quasi-tree by connecting the barycentric centers of the leaves 
and branches according to the above rules. In our example we get

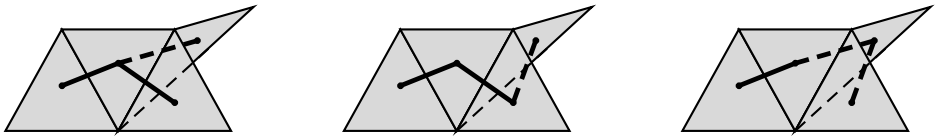
\begin{center}
\psset{unit=1.5cm}
\begin{pspicture}(0,0)(9,1.8)
\pspolygon[style=fyp, fillcolor=medium](0.5,0.5)(1.5,0.5)(1,1.4)
 \pspolygon[style=fyp, fillcolor=medium](1.5,0.5)(1,1.4)(2,1.4)
 \pspolygon[style=fyp, fillcolor=medium](1.5,0.5)(2,1.4)(2.7,1.6)
  \pspolygon[style=fyp, fillcolor=medium](1.5,0.5)(2,1.4)(2.5,0.5)
 \psline[linestyle=dashed](1.5,0.5)(2.3,1.23)
 \rput(1,0.9){\vertex}
 \rput(1.5,1.1){\vertex}
 \rput(2.2,1.3){\vertex}
 \rput(2,0.75){\vertex} 
 \psline[linewidth=1.8pt](1,0.9)(1.5,1.1)
 \psline[linewidth=1.8pt,linestyle=dashed](1.5,1.1)(2.2,1.3)
 \psline[linewidth=1.8pt](1.5,1.1)(2,0.75) 
 
 \pspolygon[style=fyp, fillcolor=medium](3.5,0.5)(4.5,0.5)(4,1.4)
 \pspolygon[style=fyp, fillcolor=medium](4.5,0.5)(4,1.4)(5,1.4)
 \pspolygon[style=fyp, fillcolor=medium](4.5,0.5)(5,1.4)(5.7,1.6)
  \pspolygon[style=fyp, fillcolor=medium](4.5,0.5)(5,1.4)(5.5,0.5)
 \psline[linestyle=dashed](4.5,0.5)(5.3,1.23)
 \rput(4,0.9){\vertex}
 \rput(4.5,1.1){\vertex}
 \rput(5.2,1.3){\vertex}
 \rput(5,0.75){\vertex} 
 \psline[linewidth=1.8pt](4,0.9)(4.5,1.1)
 \psline[linewidth=1.8pt](4.5,1.1)(5,0.75)
 \psline[linewidth=1.8pt,linestyle=dashed](5.2,1.3)(5,0.75)
 
 \pspolygon[style=fyp, fillcolor=medium](6.5,0.5)(7.5,0.5)(7,1.4)
 \pspolygon[style=fyp, fillcolor=medium](7.5,0.5)(7,1.4)(8,1.4)
 \pspolygon[style=fyp, fillcolor=medium](7.5,0.5)(8,1.4)(8.7,1.6)
  \pspolygon[style=fyp, fillcolor=medium](7.5,0.5)(8,1.4)(8.5,0.5)
 \psline[linestyle=dashed](7.5,0.5)(8.3,1.23)
 \rput(7,0.9){\vertex}
 \rput(7.5,1.1){\vertex}
 \rput(8.2,1.3){\vertex}
 \rput(8,0.75){\vertex} 
 \psline[linewidth=1.8pt](7,0.9)(7.5,1.1)
 \psline[linewidth=1.8pt,linestyle=dashed](7.5,1.1)(8.2,1.3)
 \psline[linewidth=1.8pt,linestyle=dashed](8.2,1.3)(8,0.75)
\end{pspicture}
\end{center}

\section{An algebraic proof of Dirac's theorem}

Let $G$ be a finite graph on $[n]$
without loops and multiple edges and $E(G)$ its edge set.
A {\em stable subset} of $G$ is a subset $F$ of $[n]$ 
such that $\{ i, j \} \in E(G)$ for all $i, j \in F$ with $i \neq j$. We write 
$\Delta(G)$ for the simplicial complex on $[n]$ whose faces are the stable subsets of $G$. 
It is clear that $G$ is the 1-skeleton of $\Delta(G)$,  and that if $\Gamma$ is a simplicial complex with 
$G=\skel_\Gamma(1)$, then $\Gamma$ is a subcomplex of $\Delta(G)$. Hence, in a certain sense, $\Delta(G)$ is the 
`largest' simplicial complex whose $1$-skeleton is $G$.

The following example demonstrates this concept: 

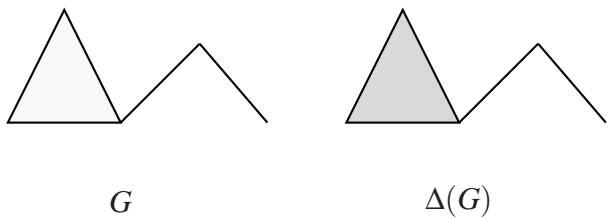
\begin{center}
\psset{unit=1.5cm}
\begin{pspicture}(0,0)(7,2.5)
\pspolygon[style=fyp](1,1)(2,1)(1.5,2)
 \psline(2,1)(2.7,1.7)
 \psline(2.7,1.7)(3.3,1)
 \rput(2,0.3){$G$}
\pspolygon[style=fyp, fillcolor=medium](4,1)(5,1)(4.5,2)
 \psline(5,1)(5.7,1.7)
 \psline(5.7,1.7)(6.3,1)
 \rput(5,0.3){$\Delta(G)$} 
\end{pspicture}
\end{center}

Recall that a graph $G$ is called {\em chordal} if each cycle of length $>3$ has a chord.  

\begin{Lemma}
\label{useful}
Let $G$ be a graph, and $\Delta$ the simplicial complex defined by $I_\Delta=I(\bar{G})$. Then 
\begin{enumerate}
\item[(a)] $\Delta=\Delta(G)$;
\item[(b)] $G=\skel_{\Delta}(1)$;
\item[(c)] $\Delta$ is a quasi-tree \iff $G$ is chordal.
\end{enumerate}
\end{Lemma}

\begin{proof}
(a) Since the 1-skeleton of $\Delta(G)=G$, it follows that $I(\bar{G})\subset I_{\Delta(G)}$. Conversely, let $F$ be a 
minimal nonface of $\Delta(G)$. If $|F|>2$, then each subset $G\subset F$
with $|G|=2$ is an edge of $G$. Therefore $F$ is a stable subset of $G$, and hence $F\in \Delta(G)$, a contradiction. 
Thus for every minimal nonface $F$ of $\Delta(G)$ one has $|F|=2$. This shows that $I_{\Delta(G)}=I(\bar{G})$. 
Therefore, $\Delta=\Delta(G)$.

(b) follows from Lemma \ref{trivial} (or from (a) and the remarks preceding this lemma).

(c) The theorem of  Fr\"oberg \cite{F} guarantees that the complementary graph 
$G$ of $\bar{G}$ is a chordal graph  if and only if $I(\bar{G})=I_{\Delta}$ has 
a $2$-linear resolution.  By  Theorem \ref{duality},  $\reg(I_{\Delta}) = \projdim I_{\Delta^\vee} + 1$, and so the 
ideal $I(\bar{G})$ has a $2$-linear
resolution if and only if $\projdim I_{\Delta^\vee} = 1$. Since by Lemma \ref{true}, $I_{\Delta^\vee} =I(\Delta^c)$, 
the assertion follows from  Corollary \ref{projdim}. 
\end{proof}

For our proof of Dirac's theorem we also need

\begin{Lemma}
\label{flaglemma}
A quasi-tree is a flag complex.
\end{Lemma}

\begin{proof}
Let $\Delta$ be a quasi-tree on $[n]$ and fix a leaf ordering
of the facets $F_1, \ldots, F_t$ of $\Delta$.
We work induction on $t$.    
Let $t > 2$.  Since $\Delta' = \langle F_1,\ldots,F_{t-1}\rangle$
is a quasi-tree, by assumption of induction 
it follows that $\Delta'$ is flag.
Let $F_k$ with $k <t$ be a branch of $F_t$. 
Then $\Delta'$ consists of all faces $G$ of $\Delta$
with $G \cap (F_t \setminus F_k) = \emptyset$.
Suppose  $H$ is a minimal nonface of $\Delta$ having at least
three elements of $[n]$.  We then show that $H$
is a minimal nonface of $\Delta'$, i.e., 
$H \cap (F_t \setminus F_k) = \emptyset$.
Since $H$ is a nonface, there is $p \in H$
with $p \not\in F_t$.  If $q \in F_t$ belongs to
$H$, then $\{ p, q \} \in \Delta$.
Thus there is $F_j$ with $j \neq t$ 
such that $\{ p, q \} \subset F_j$.
Hence $\{ q \} \subset F_t \cap F_j$.
Thus $q \in F_k$.  
Hence $H \cap (F_t \setminus F_k) = \emptyset$,
as desired.
\end{proof}

\begin{Theorem}[Dirac]
\label{Dirac}
A finite graph $G$ on $[n]$ is a chordal graph if and only if
$G$ is the $1$-skeleton of a quasi-tree on $[n]$.
\end{Theorem} 

\begin{proof}
The statements (b) and (c) of Lemma \ref{useful} imply that a chordal graph is the 1-skeleton of quasi-tree. 
Conversely, suppose that $G$ is the 1-skeleton of a quasi-tree $\Gamma$. 
Since by Lemma \ref{flaglemma}, $\Gamma$ is flag, the ideal $I_\Gamma$ is generated by all monomials $x_F$ with 
$|F|=2$ and $F\not\in \Gamma$. This shows that $I_\Gamma=I(\bar{G})$, and so $\Gamma= 
\Delta(G)$, by Lemma \ref{useful}(a).  Hence $G$ is chordal by Lemma \ref{useful}(c).
\end{proof}

\begin{Corollary}
\label{hibifound}
Let $\Delta$ be a quasi-tree, and $F$ a leaf of $\Delta$. Then $\langle {\mathcal F}(\Delta)\setminus \{F\}\rangle$ is 
again a quasi-tree.
\end{Corollary}

\begin{proof}
Let $\Delta' = \langle {\mathcal F}(\Delta)\setminus \{F\}\rangle$.  Let $G$ be the $1$-skeleton of $\Delta$ and $G'$ 
the $1$-skeleton of $\Delta'$.  Then $G'$ is obtained by removing all free vertices of $F$ and all edges containing 
these vertices from $G$.  Since $G$ is chordal by Theorem \ref{Dirac}, it follows that $G'$ is also chordal.  Hence 
again by Theorem \ref{Dirac}, $\Delta'$ is a quasi-tree.
\end{proof}

We conclude this section with a  sort of higher Dirac theorem.

\begin{Theorem}
\label{higher}
Let $\Delta$ be a pure $\ell$-dimensional  simplicial complex on the vertex set $[n]$, and $\Gamma$ its $1$-skeleton. 
Then the following conditions are equivalent:
\begin{enumerate}
\item[(a)] $\Delta$ is the $\ell$-skeleton of a quasi-tree;
\item[(b)] 
\begin{enumerate} 
\item[(i)] $\Gamma$ is a chordal graph;
\item[(ii)] $\Delta$ is the $\ell$-skeleton of $\Delta(\Gamma)$.
\end{enumerate}
\end{enumerate}
\end{Theorem}

\begin{proof}
(b)\implies (a) follows from Lemma \ref{useful}(c). For the implication (a)\implies (b), suppose that $\Delta$ is the 
$\ell$-skeleton of the quasi-tree $\Sigma$. Then $\Gamma$ is also the $1$-skeleton of  $\Sigma$. As in the proof of 
Theorem \ref{Dirac} we conclude  that $\Sigma=\Delta(\Gamma)$. This implies (b)(ii). Finally, by Dirac's theorem 
$\Gamma$ is chordal.
\end{proof}

\section{Powers of facet ideals related to quasi-trees}

We now consider powers of facet ideals of complementary simplicial complexes of skeletons of quasi-trees. We first 
show that such ideals have linear quotients.

\begin{Theorem}
\label{linear quotient}
Let $\Delta$ be a quasi-tree of dimension $d-1$.  Then 
$I=I(\overline{\skel_\Delta(\ell)})$ has linear quotients for any $\ell\leq d-1$. In particular, $I$ has a linear 
resolution.
\end{Theorem}

\begin{proof}
Let $I_\Gamma=I$ and  $I_{\Gamma'}=I(\overline{\skel_\Delta(1)})$. Since by the Lemma  \ref{flaglemma}, $\Delta$ is 
flag  we have $I_\Delta=I_{\Gamma'}$. In  \cite{HHZ} we showed using Dirac's theorem  that  $I_{\Gamma'}$ has linear 
quotients. 
By Corollary \ref{application}, $I$ has linear quotients, too. 
\end{proof}

In \cite{HHZ} a certain converse of Theorem \ref{linear quotient} is shown for $\ell=1$, namely, that if  $I$ is a 
monomial ideal generated in degree 2 and has linear quotients, then there exists a quasi-tree $\Delta$ such that $I= 
I(\overline{\skel_\Delta(1)})$.  However, for $\ell>1$, such a  converse is not true:  let 
$\Delta=\langle\{1,2,3\},\{3,4,5\},\{2,4,6\}\rangle$, and  $I=I(\bar{\Delta})$. Then $I$ has linear quotients. 
However, if $I=I(\overline{\skel_\Gamma(2))}$, then  $\Delta=\skel_\Gamma(2)$. In particular, $\dim \Gamma\geq 2$. If 
$\dim \Gamma>2$, then $\skel_\Gamma(2)$ contains at least 4 facets. But $\Delta$ has only 3 facets. Thus 
$\dim\Gamma=2$, and hence $\Gamma=\Delta$. But $\Delta$ is not a quasi-tree.

The main theorem of this section is the following

\begin{Theorem}
\label{wehope}
Let $\Delta$ be a quasi-tree of dimension $d-1$. Then for any $\ell\leq d-1$, all powers of 
$I=I(\overline{\skel_\Delta(\ell)})$ have a linear resolution. 
\end{Theorem}

We need the following two lemmata.

\begin{Lemma}
\label{combinatorialfeeling}
Let $2 \leq \ell \leq d-1$,  
$I_1=I(\overline{\skel_\Delta(1)})$ and
$I_\ell=I(\overline{\skel_\Delta(\ell)})$.
Then $I_\ell$ is generated by all squarefree monomials $u$ of degree $\ell +1$ such that $u$ is divided by a monomial 
generator of $I_1$.
\end{Lemma}

\begin{proof}
Let $u = x_F$ be a squarefree monomial of degree $\ell + 1$.
If $u$ is divided by a monomial generator $x_{i}x_{j}$ of $I_1$,
then $F$ contains the $2$-element subset $\{i, j\} \not\in \Delta$.
Thus $F \not\in \Delta$ and $u$ is a monomial generator of $I_\ell$.
Conversely, suppose that $u$ is divided by no monomial generator of $I_1$.
Then each $2$-element subset of $F$ is a face of $\Delta$.  
Since $\Delta$ is flag, it follows that $F$ is a face of $\Delta$.
Thus $u \not\in I_\ell$. 
\end{proof}

Given integer vectors $a = (a_1, \ldots, a_n)$ and $b = (b_1, \ldots, b_n)$, we write $a \leq b$ if $a_i \leq b_i$ for 
all $i$.  
Let $I\subset S$ be an arbitrary monomial ideal, and $a=(a_1,\ldots, a_n)$ an integer vector with each $a_i\geq 0$. We 
write $I^{\leq a}$ for the monomial ideal generated by all $u=x^b \in G(I)$ with $b \leq a$.  Here $x^b = 
x_1^{b_1}\cdots x_n^{b_n}$ if $b = (b_1, \ldots, b_n)$.

\begin{Lemma}
\label{algebraicfeeling}
Let $I\subset S$ be a monomial ideal,  
\[ 
\Bbb F\: 0\To F_p\To   F_{p-1}\To\cdots\To F_1\To F_0\To S/I\To 0
\]
the multigraded minimal free resolution of $I$ with $F_i=\Dirsum_jS(-q_{ij})$, and $\Bbb G$ the  subcomplex of $\Bbb 
F$ with 
\[
G_i=\Dirsum_{q_{ij}\leq a}S(-q_{ij}).
\]
Then $\Bbb G$ is a  multigraded minimal free resolution of $I^{\leq a}$.  In particular, if $I$ has a linear 
resolution, then so does $I^{\leq a}$. 
\end{Lemma}

\begin{proof}
It is clear that  $H_0(\Bbb G)=S/I^{\leq a}$. Thus it remains to show that $\Bbb G$ is acyclic. 
We proceed by induction on the  homological degree.  Suppose that our claim is true up to homological degree $i$, and 
let $r$ be a  multihomogeneous element belonging to a minimal set of generator  of the kernel of  $\Bbb G_i\to \Bbb 
G_{i-1}$. Let $v$ be the multidegree of $r$. It is known \cite{BH1} that $v\leq a$.

Now $r$ belongs to the kernel $C$ of $\Bbb F_i\to \Bbb F_{i-1}$ as well.  Let $\{c_1,\ldots, c_m\}$ be the minimal set 
of generators of $C$ corresponding to the chosen basis of $F_{i+1}$. Then  $r=\sum_i h_ic_i$  where each $h_ic_i$ has 
the same multidegree  as $r$.  It is then  clear $h_i\neq 0$ only if the multidegree of $c_i$ is bounded by $a$. Hence 
$r$ belongs to the image of $G_{i+1}\to G_i$, as required.
\end{proof}

\begin{proof}[Proof of Theorem \ref{wehope}]
Let $I=I(\overline{\skel_\Delta(1)})$ and $J=I(\overline{\skel_\Delta(\ell)})$. 
By Lemma \ref{combinatorialfeeling} it follows that $J= (I_{\langle \ell+1\rangle})^{\leq (1,\ldots,1)}$, where for 
some graded ideal $L$, we denote by $L_{\langle j\rangle}$ the ideal generated by the elements of the $j$th graded 
component of $L$. Note that $J^k=((I^k)_{\langle k(\ell+1)\rangle})^{(k,\ldots,k)}$. By \cite[Theorem 3.2]{HHZ}, $I^k$ 
has a linear resolution. Hence $(I^k)_{\langle k(\ell+1)\rangle}$ has a linear resolution. Then Lemma 
\ref{algebraicfeeling} guarantees that $J^k$ has a linear resolution.

\end{proof}


\begin{thebibliography}{99}

\bibitem{BH} W.\ Bruns and J.\ Herzog, ``Cohen--Macaulay
rings,''  Revised Edition, Cambridge University Press, 
Cambridge, 1996. 


\bibitem{BH1} W.\ Bruns and J.\ Herzog, On multigraded resolutions, {\it Math.\ Proc.\ Camb.\ Phil.\ Soc.} {\bf 118} 
(1995), 245 -- 257.


\bibitem{D} G.\ A.\ Dirac, On rigid circuit graphs, {\it Abh.\ Math.\ Sem.\ Univ.\ Hamburg}, {\bf 38} (1961), 71 -- 
76. 

\bibitem{Fa} S.\ Faridi, The facet ideal of a simplicial complex, {\it Manuscripta Mat.} {\bf 109} (2002), 159 -- 174.   
 

\bibitem{ER} J.\ A.\ Eagon and V.\ Reiner, Resolutions of Stanley--Reisner rings and Alexander duality, {\it J.\ Pure 
and Appl.\ Alg.} {\bf 130} (1998), 265 -- 275. 

\bibitem{F} R.\ Fr\"oberg,  On Stanley--Reisner rings,  in: Topics in algebra, {\it Banach Center Publications}, {\bf 
26} (2), (1990), 57 -- 70.

\bibitem{HHZ} J.\ Herzog, T.\ Hibi and X.\ Zheng, Monomial ideals whose powers have a linear resolution, Preprint 
2002. 

\bibitem{T} N.\ Terai, Generalization of Eagon--Reiner theorem and $h$-vectors of graded rings, Preprint 2000.


\bibitem{Z} X.\ Zheng,  Resolutions of facet ideals, to appear in Comm.\ Alg. 

\end{thebibliography}
\end{document}